\def\SScurfname{ \jobname.tex}
\def\SSone{}
\everypar={\SSone}

\documentclass[11pt, twoside]{article}
\usepackage{bbm}
\usepackage{amssymb}

\title{\bf A probabilistic approach to enumeration of Gessel walks}

\author{Sun Ping
\footnote{\it E-mail address: plsun@mail.neu.edu.cn.}}

\date{}

\begin{document}

\maketitle

\centerline{{\footnotesize \it Department of Mathematics,
Northeastern University, Shenyang, 110004, P R China}}

\vspace{10pt}

{\bf Abstract.} We consider Gessel walks in the plane starting at
the origin $(0, 0)$ remaining in the first quadrant $i, j \geq 0$
and made of West, North-East, East and South-West steps. Let $F(m;
n_1, n_2)$ denote the number of these walks with exact $m$ steps
ending at the point $(n_1, n_2)$, Petkov\v{s}ek and Wilf posed
several analogous conjectures similar to the famous Gessel's
conjecture.

We establish a probabilistic model of Gessel walks which is
concerned with the problem of vicious walkers. This model helps us
to obtain the linear homogeneous recurrence relations with binomial
coefficients for both $F(n+k+r;n+k-r,n)$ and $F(n+2k; n, 0)$.
Precisely, $\frac{n! k! (n+k+1)!}{(2n+2)!} F(2n+2k;0,n)$ is a
polynomial with all integer coefficients which leading term is
$2^{3k-2} n^{2k-2}$, and $\frac{k! (k+1)!}{n+1} F(n+2k;n,0)$ is a
polynomial with all integer coefficients which leading term is
$n^{2k-1}$. Hence two conjectures of Petkov\v{s}ek and Wilf are
solved.

\vspace{10pt}


\vspace{5pt}


\section{Introduction}

\paragraph*{}

Let $F(m; n_1, n_2)$ be the number of lattice walks (called Gessel
walks in literatures) from the origin $(0,0)$ to $(n_1, n_2)$,
consisting of $m$ steps that can be West, North-East, East and
South-West (abbreviated, W, N-E, E, and S-W), remaining in the first
quadrant $\{(i, j); i \geq 0, j \geq 0 \}$, the obvious recurrence
form is
$$
F(m; n_1, n_2) = F(m-1; n_1+1, n_2) + F(m-1; n_1-1, n_2) + F(m-1;
n_1+1, n_2+1)
$$
$$  + F(m-1; n_1-1, n_2-1), \hspace{5pt} m \geq 1, n_1 \geq 0, n_2 \geq 0,  \eqno(1)
$$
with the initial condition
$$
F(0;n_1,n_2) = \cases{ 1, & $n_1 = n_2 = 0$, \cr
                         0, & otherwise, } \eqno(2)
$$
and the boundary condition
$$
F(m; n_1, n_2) = 0, \hspace{2pt} n_1 < 0 \hspace{5pt} or
\hspace{5pt} n_2 < 0. \eqno(3)
$$

It's clear that $F(m;n_1,n_2) \neq 0$ only if \cite{PW}
$$
m \equiv n_1(mod \hspace{3pt} 2), \hspace{5pt}   n_1 \leq m,
\hspace{5pt} n_2 \leq \frac{1}{2}(m+n_1).
$$

\paragraph*{}

The enumeration of Gessel walks is a difficult problem in
combinatorics, since Ira Gessel conjectured empirically that
$F(2n;0,0)$ has a closed form in 2001 \cite{KZ}:
$$
F(2n;0,0) = 16^n \frac{(5/6)_n (1/2)_n}{(2)_n (5/3)_n}, \eqno(4)
$$
where $(a)_n = a(a+1) \dots (a+n-1)$ is the Pochhammer symbol,
recently A.Bostan and M.Kauers \cite{BK} have shown that the
generating function of $F(m;n_1n_2)$ is algebraic, then is
holonomic, that is $F(m;n_1,n_2)$ satisfies linear recurrences with
polynomial coefficients, M.Bousquet-Melou and M.Mishna \cite{ME}
described a complete classification of all the 256 step sets $S
\subset  \{ -1,0,1 \}^2 \backslash \{0,0\}$, they can solve all the
problems with a unified way except Gessel walks. It is remarkable
that M.Kauers, C.Koutschan and D.Zeilberger \cite{KK} gave a
computer-aided proof of the Gessel conjecture in 2008, their proof
is accomplished by computing a homogeneous linear recurrence in $m$
for $F(m;0,0)$, however the recurrence has order $32$, polynomial
coefficients of degree $172$ and involves integers with up to $385$
decimal digits, and they believe it is very possible that a short
human proof of Gessel conjecture does not exist. Furthermore,
M.Petkov\v{s}ek and H.S.Wilf \cite{PW} derived a functional equation
of the generating function of $F(m;n_1,n_2)$ and turned it into an
infinite lower-triangular system of linear equations satisfied by
the values of $F(m;n_1,0)$ and $F(m;0,n_2) + F(m;0,n_2 -1)$ using
the kernel method, and these values are expressed as determinants of
lower-Hessenberg matrices with unit super diagonals whose non-zero
entries are products of two binomial coefficients, meanwhile they
posed some similar conjectures {\bf C1 - C4} as follows:

\vspace{10pt}

{\bf C1}. {\it For} $n \geq 1$,
$$
F(2n;0,1)=16^n \frac{(1/2)_n}{(3)_n} \left( \frac{5}{27}
\frac{(7/6)_n}{(7/3)_n} + \frac{(111n^2 + 183n - 50)}{270}
\frac{(5/6)_n}{(8/3)_n} \right), \eqno(5.1)
$$
{\it and more generally, for} $n \geq k$,
$$
F(2n;0,k)=16^n \frac{(1/2)_n}{(k+2)_n} \left(
\frac{(7/6)_n}{((3k+4)/3)_n} p_k(n) + \frac{(5/6)_n}{((3k+5)/3)_n}
q_k(n) \right), \eqno(5.2)
$$
{\it where $p_k(n)$ is a polynomial of degree $2k-2$ and $q_k(n)$ is
a polynomial of degree $2k$ }.

\vspace{5pt}

{\bf C2}. {\it For $n \geq 0$},
$$
\left. \begin{array}{l} F(2n;0,n)=4^n \frac{(3/2)_n}{(2n+1)(2)_n} =
\frac{4^n (1/2)_n}{(2)_n}, \\
F(2n+2;0,n)=\frac{2^{2n+1}(n+1)(3/2)_n}{(3)_n}, \\
F(2n+4;0,n)=\frac{4^n(n+1)(8n^2 + 32n + 33)(3/2)_n}{3(4)_n}, \\
F(2n+6;0,n)=\frac{4^{n-1}(n+1)(64n^4 + 672n^3 + 2648n^2 + 4641n +
3060)(3/2)_n}{9(5)_n},
\end{array} \right \} \eqno(6.1)
$$
{\it and then}
$$
F(2n+2k;0,n) = 4^n \frac{(3/2)_n}{(k+2)_n}r_k(n), \eqno(6.2)
$$
{\it where $r_k(n)$ is a polynomial of degree $2k-1$ divisible by
$n+1$ for $k \geq 1$ and $r_0(n)=1/(2n+1).$}

\vspace{5pt}

{\bf C3}. {\it Let $g(n)= F(2n+1;1,0)$, it satisfies the second
order recurrence}
$$
(n+3)(3n+7)(3n+8)\hspace{2pt} g(n+1)
$$
$$
\hspace{-2cm} - 8 (2n+3)(18n^2 + 54n+35)\hspace{2pt} g(n) \eqno(7)
$$
$$ + 256 n(3n+1)(3n+2)\hspace{2pt} g(n-1)=0.
$$

\vspace{5pt}

{\bf C4}. {\it For $n \geq 0$},
$$
\left. \begin{array}{l} F(n;n,0)=1, \\
F(n+2;n,0)=\frac{1}{2}(n+1)(n+4), \\
F(n+4;n,0)=\frac{1}{12}(n+1)(n^3 + 15n^2 + 74n + 132), \\
F(n+6;n,0)=\frac{1}{144}(n+1)(n^5 + 32n^4 + 407n^3 + 2620n^2 + 8604n
+ 12240),
\end{array} \right\} \eqno(8.1)
$$
{\it and then}
$$
F(n+2k;n,0) = s_k(n), \eqno(8.2)
$$
{\it where $s_k(n)$ is a polynomial of degree $2k$ with leading
coefficient $\frac{1}{k! (k+1)!}$, which is divisible by $n+1$ for
$k \geq 1$.}

\paragraph*{}

In \cite{KK} Kauers, Koutschan and Zeilberger have given  a computer
proof of the equality (5.1) in {\bf C1} and {\bf C3}, they also
found a recurrence for $F(2n;2,0)$. By the method of generation
Petkov\v{s}ek and Wilf also obtained the following result \cite{PW}:
$$
F(n_1;n_1,n_2)={n_1 \choose n_2}, \hspace{5pt} n_1 \geq n_2,
\eqno(9)
$$
$$
F(2n_2 - n_1;n_1,n_2)=\frac{n_1 + 1}{2n_2 - n_1 +1}{2n_2 - n_1 +1
\choose n_2 + 1}, \hspace{5pt} n_1 \leq n_2. \eqno(10)
$$

\paragraph*{}

In Section $2$ we will introduce a new probabilistic method for
lattice path enumeration which is different from that used by Ira
Gessel in \cite{GE}, we make connections to Gessel walks and to
pairs of non-crossing lattice paths, or vicious walkers. In Section
$3$ the conjectures {\bf C2} and {\bf C4} are proved by using the
probabilistic methods, the techniques that we employ are basic.
Finally, we discuss related problems on the Gessel conjecture $(4)$
in Section $4$. Gessel number $F(2n;0,0)$ is shown to be the number
of 2-watermelons of length $2n$ with half-wall, also the number of
pairs of non-crossing Dyck paths and free Dyck paths of length $2n$.
It seems it is possible to give a human proof of Gessel conjecture
by applying the theory of vicious walkers.

\newpage

\section{A probabilistic approach to Gessel walks}

\paragraph*{}

Noticed the four directions of Gessel walks, suppose a random
particle from $(x_0, y_0)$ takes one of W, N-E, E and S-W steps with
probability $1/4$ to the point $(x_0+x,y_0+y)$, then $(x, y)$ is  a
random vector with the following distribution:

\vspace{5pt}

\begin{table}[h]
\begin{flushright}
\begin{tabular}{c|cccr}
{\it x $\backslash$ y}   &   -1  &  0     &  1    &  \\ \cline{1-4}
 -1                      &   1/4 &  1/4   &  0    & \hspace{4.5cm} (11) \\
 1                       &    0  &  1/4   &  1/4  &
\end{tabular}
\end{flushright}
\end{table}

\vspace{5pt}

Denote $(x_1,y_1), (x_2, y_2), \cdots$ be independent and
identically distributed (abbreviated i.i.d) random vectors with the
distribution $(11)$, the partial sums be
$$
X_j = \sum_{i=1}^j x_i, \hspace{10pt} Y_j = \sum_{i=1}^j y_i,
\hspace{10pt} X_0 \equiv 0, Y_0 \equiv 0,
$$
then we have
$$
F(m;n_1,n_2) = 4^m P \{ X_j \geq 0, Y_j \geq 0, 1 \leq j \leq m, X_m
= n_1, Y_m = n_2 \}. \eqno(12)
$$

It's inconvenient to enumerate Gessel walks using the representation
$(12)$ because $x, y$ are dependent, in this paper we introduce
i.i.d random variables $\xi, \eta$ with Bernoulli distribution:
$$
P\{\xi = -1\} = P\{\xi = 1\} = P\{\eta = -1\} = P\{\eta = 1\} =
\frac{1}{2}, \eqno(13)
$$
and $(\xi_1,\eta_1), (\xi_2, \eta_2), \cdots $ are i.i.d random
vectors with the distribution $(13)$, the partial sums are
$$
S_j = \sum_{i=1}^j \xi_i, \hspace{10pt} T_j = \sum_{i=1}^j \eta_i,
\hspace{10pt} S_0 \equiv 0, T_0 \equiv 0,
$$
it is interesting that the distribution of $(\xi, \frac{\xi -
\eta}{2} )$ is identical with $(11)$ of $(x,y)$, then $(12)$ implies

\vspace{5pt}

{\bf Lemma 1}. {\it For $m \geq 0, n_1 \geq 0 , n_2 \geq 0$}
$$
F(m;n_1,n_2) = 4^m P \{ S_j \geq 0, S_j \geq T_j, 1 \leq j \leq m,
S_m = n_1, T_m = n_1 - 2n_2 \}, \eqno(14)
$$
{\it specially, Gessel numbers have the following representation}
$$
F(2n;0,0)=4^{2n} P \{ S_j \geq 0, S_j \geq T_j, 1 \leq j \leq 2n,
S_{2n} = T_{2n} = 0 \}. \eqno(15)
$$

\paragraph*{}

We note that the equalities $(9), (10)$ of Petkov\v{s}ek and Wilf
can be derived immediately by using the probabilistic representation
of $F(m;n_1,n_2)$. Because $S_{n_1} = n_1$ means that
$\xi_1=\cdots=\xi_{n_1}=1$, and $T_{n_1} = n_1 - 2n_2$ implies that
exactly $n_1 - n_2 (n_1 \geq n_2)$ of $\{ \eta_1=1 \}, \cdots, \{
\eta_{n_1}=1 \} $ occur, from $(14)$ we have
$$
F(n_1;n_1,n_2)=2^{n_1} P \{ T_{n_1}=n_1 -2n_2 \} = {n_1 \choose
n_2}.
$$

On the other hand, $T_{2n_2 - n_1}=n_1-2n_2$ means $\eta_1 = \cdots
= \eta_{2n_2-n_1}=-1$, and $S_{2n_2-n_1}=n_1$ implies that exactly
$n_2-n_1 (n_1 \leq n_2)$ of $\{ \xi_1=-1 \}, \cdots, \{
\xi_{2n_2-n_1} = -1 \}$ occur, from the classic result of W. Feller
(\cite{FE}, Ch3) we have
$$
F(2n_2-n_1;n_1,n_2)=2^{2n_2-n_1} P \{ S_j \geq 0, 1 \leq j \leq
2n_2-n_1, S_{2n_2-n_1}=n_1 \}
$$
$$
= \frac{n_1 + 1}{2n_2 - n_1 +1}{2n_2 - n_1 +1 \choose n_2 + 1}.
$$

\paragraph*{}

It's interesting that Lemma $1$ shows $F(m;n_1,n_2)$ is the number
of pairs of non-crossing lattice paths  $(S,T)$ of length $m$
starting from $(0,0)$ with steps of the form North-East and
South-East, where $S$ running above (may touching) the $x$-axis and
ending at $(m,n_1)$, $T$ running below (may touching) $S$ ending at
$(m,n_1-2n_2)$. Figure $1$ is an illustration. Moreover, in Section
$4$, we will show that Gessel walks is the special case of two
vicious walkers with half-wall.

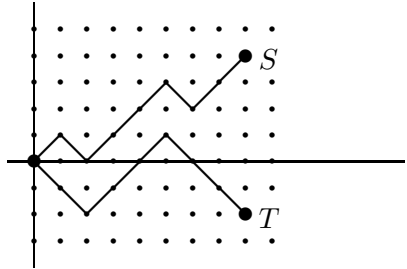
\begin{figure}[ht]
\begin{center}
\begin{picture}(100,100)
\put(10,0){\line(0,1){100}} \put(0,40){\line(1,0){150}}
\put(10,40){\circle*{5}} \put(20,40){\circle*{2}}
\put(30,40){\circle*{2}} \put(40,40){\circle*{2}}
\put(50,40){\circle*{2}} \put(60,40){\circle*{2}}
\put(70,40){\circle*{2}} \put(80,40){\circle*{2}}
\put(90,40){\circle*{2}} \put(100,40){\circle*{2}}
\put(10,30){\circle*{2}} \put(20,30){\circle*{2}}
\put(30,30){\circle*{2}} \put(40,30){\circle*{2}}
\put(50,30){\circle*{2}} \put(60,30){\circle*{2}}
\put(70,30){\circle*{2}} \put(80,30){\circle*{2}}
\put(90,30){\circle*{2}} \put(100,30){\circle*{2}}
\put(10,20){\circle*{2}} \put(20,20){\circle*{2}}
\put(30,20){\circle*{2}} \put(40,20){\circle*{2}}
\put(50,20){\circle*{2}} \put(60,20){\circle*{2}}
\put(70,20){\circle*{2}} \put(80,20){\circle*{2}}
\put(90,20){\circle*{5}} \put(10,10){\circle*{2}}
\put(20,10){\circle*{2}} \put(30,10){\circle*{2}}
\put(40,10){\circle*{2}} \put(50,10){\circle*{2}}
\put(60,10){\circle*{2}} \put(70,10){\circle*{2}}
\put(80,10){\circle*{2}} \put(90,10){\circle*{2}}
\put(100,10){\circle*{2}} \put(10,50){\circle*{2}}
\put(20,50){\circle*{2}} \put(30,50){\circle*{2}}
\put(40,50){\circle*{2}} \put(50,50){\circle*{2}}
\put(60,50){\circle*{2}} \put(70,50){\circle*{2}}
\put(80,50){\circle*{2}} \put(90,50){\circle*{2}}
\put(100,50){\circle*{2}} \put(10,60){\circle*{2}}
\put(20,60){\circle*{2}} \put(30,60){\circle*{2}}
\put(40,60){\circle*{2}} \put(50,60){\circle*{2}}
\put(60,60){\circle*{2}} \put(70,60){\circle*{2}}
\put(80,60){\circle*{2}} \put(90,60){\circle*{2}}
\put(100,60){\circle*{2}} \put(10,70){\circle*{2}}
\put(20,70){\circle*{2}} \put(30,70){\circle*{2}}
\put(40,70){\circle*{2}} \put(50,70){\circle*{2}}
\put(60,70){\circle*{2}} \put(70,70){\circle*{2}}
\put(80,70){\circle*{2}} \put(90,70){\circle*{2}}
\put(100,70){\circle*{2}} \put(10,80){\circle*{2}}
\put(20,80){\circle*{2}} \put(30,80){\circle*{2}}
\put(40,80){\circle*{2}} \put(50,80){\circle*{2}}
\put(60,80){\circle*{2}} \put(70,80){\circle*{2}}
\put(80,80){\circle*{2}} \put(90,80){\circle*{5}}
\put(10,90){\circle*{2}} \put(20,90){\circle*{2}}
\put(30,90){\circle*{2}} \put(40,90){\circle*{2}}
\put(50,90){\circle*{2}} \put(60,90){\circle*{2}}
\put(70,90){\circle*{2}} \put(80,90){\circle*{2}}
\put(90,90){\circle*{2}} \put(100,90){\circle*{2}}

\thicklines \put(10,40){\line(1,1){10}} \put(20,50){\line(1,-1){10}}
\put(30,40){\line(1,1){30}} \put(60,70){\line(1,-1){10}}
\put(70,60){\line(1,1){20}} \put(10,40){\line(1,-1){20}}
\put(30,20){\line(1,1){30}} \put(60,50){\line(1,-1){30}}
\put(95,75){$S$} \put(95,15){$T$}

\end{picture}
\caption{$F(8;4,3)=\# \{(S,T): S \geq 0, S \geq T, S_8=4,T_8=-2 \}$}
\end{center}
\end{figure}

\paragraph*{}

Now we consider the case of Gessel walks confined in the half-plane.
Let $G(m;n_1,n_2)$ be the number of lattice walks from $(0,0)$ to
$(n_1, n_2)$, consisting of $m$ steps that can be W, N-E, E and S-W,
remaining in the first and second quadrant $\{(i, j); j \geq 0 \}$,
from lemma 1 it follows that
$$
G(m,n_1,n_2) = 4^m P \{ S_j \geq T_j, 1 \leq j \leq m, S_m = n_1,
T_m = n_1 - 2n_2 \},
$$
where $S_j, T_j$ are defined as above.

The closed formula for $G(m;n_1n_2)$ follows immediately from a
result of stars with fixed end points and without a wall on the
problem of vicious walkers (\cite{KG}, Theorem $1$), it can also be
derived by using Lindstr\"{o}m-Gessel-Viennot theorem of
non-intersecting lattice paths (\cite{GV}, Corollary $2$. Also
\cite{KG}, Proposition A2). In fact, using reflection principle it
is not difficult to evaluate $G(m;n_1,n_2)$ directly as follows.
Since the number of lattice walks from $(0,0)$ to $(n_1, n_2)$,
consisting of $m$ steps that can be W, N-E, E and S-W, is equal to
$4^m P \{ S_m = n_1, T_m = n_1 - 2n_2 \}$, thus the independent of
$\xi$ and $\eta$ implies it is equal to
$$
2^m P \{ S_m = n_1 \} \times 2^m P \{ T_m = n_1 - 2n_2 \} = {m
\choose \frac{m-n_1}{2}} {m \choose \frac{m-n_1}{2} + n_2}.
$$

\begin{figure}[ht]
\begin{center}
\begin{picture}(100,100)

\thinlines \put(0,40){\vector(1,0){150}} \put(30,40){\circle*{5}}
\put(10,20){\circle*{5}} \put(50,30){\circle{4}}
\put(80,80){\circle*{5}} \put(90,80){$(n_1,n_2)$}
\put(10,30){\line(1,0){130}} \put(110,20){$y=-1$}
\put(10,50){$(0,0)$} \put(10,50){$(0,0)$} \put(-35,25){$(-2,-2)$}

\thicklines \put(30,40){\line(1,1){10}} \put(40,50){\line(1,0){10}}
\put(50,50){\line(-1,-1){10}} \put(40,40){\line(1,0){20}}
\put(60,40){\line(-1,-1){10}}

\put(50,30){\line(1,0){20}} \put(70,30){\line(1,1){20}}
\put(95,55){\circle*{1.3}} \put(100,60){\circle*{1.3}}
\put(105,65){\circle*{1.3}}

\put(10,20){\line(-1,-1){10}} \put(0,10){\line(1,0){10}}
\put(10,10){\line(1,1){10}} \put(20,20){\line(1,0){20}}
\put(40,20){\line(1,1){10}}

\end{picture}
\caption{A bijection}
\end{center}
\end{figure}
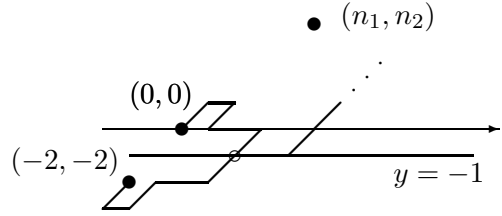

On the other hand, by Andr\'{e}'s reflection principle, there is a
bijection between walks from $(0,0)$ to $(n_1,n_2)$ which cross the
$x$-axis and walks from $(-2,-2)$ to $(n_1,n_2)$, see Figure $2$ for
an example. It is obvious the number of lattice walks from $(-2,-2)$
to $(n_1,n_2)$, consisting of $m$ steps that can be W, N-E, E and
S-W, is equal to
$$
{m \choose \frac{m-n_1}{2}-1} {m \choose \frac{m-n_1}{2} + n_2 +1},
$$
therefore, we have the following result:

\vspace{15pt}

{\bf Lemma 2}. {\it For $m \geq 0, n_2 \geq 0$}
$$
G(m;n_1,n_2) = 4^m P \{ S_j \geq T_j, 1 \leq j \leq m, S_m = n_1,
T_m = n_1 - 2n_2 \}
$$
$$
= \frac{n_2 + 1}{m+1} {m + 1 \choose \frac{m-n_1}{2}} {m + 1 \choose
\frac{m-n_1}{2} + n_2 +1}. \eqno(16)
$$

\section{Proof of the conjectures of Petkov\v{s}ek and Wilf}

\paragraph*{}

In this section, we present a probabilistic approach to the problems
of Petkov\v{s}ek and Wilf. That is, we apply Lemma $1$ to establish
the recurrences with binomial coefficients of $F(n+k+r;n+k-r,n)$ and
$F(n + 2k;n,0)$, the main results of this paper are stated as
follows:

\vspace{10pt}

{\bf Theorem 3.} {\it For $n \geq 0, \hspace{2pt} k \geq 1,
\hspace{5pt} r \leq n+k$},
$$
F(n+k+r;n+k-r,n) = \sum_{t=0}^{k-2} \left[ {n-k+r+2 \choose n - t} -
{n-k+r+2 \choose n+t+3} \right] F(2k-2;2t+2,t)
$$
$$
+ \sum_{s=0}^{k-1} \sum_{t=0}^{k-1}{n-k+r+2 \choose s+1} \left[
{n-k+r+2 \choose n-t+1} - {n-k+r+2 \choose n+t+2} \right]
F(2k-2;2t,s+t). \eqno(17)
$$

\vspace{15pt}

{\bf Corollary 4.}

\centerline{\it The conjecture {\bf C2} of Petkov\v{s}ek and Wilf is
true.}

\vspace{15pt}

{\bf Theorem 5.} {\it For $n \geq 0, k \geq 1$},
$$
F(n + 2k;n,0) = \sum_{s=1}^k \sum_{t=0}^s  \frac{s-t+1}{n+2} {n+2
\choose s+1}  {n+2 \choose t} F(2k-1;2s-1,s-t). \eqno(18)
$$

\vspace{10pt}

{\bf Corollary 6.}

\centerline{\it The conjecture {\bf C4} of Petkov\v{s}ek and Wilf is
true.}

\vspace{20pt}

{\bf Proof of Theorem 3.} We start by observing that Lemma $1$
implies
$$
F(n+k+r;n+k-r,n)=4^{n+k+r}P\{S_j \geq 0, S_j \geq T_j, S_{n+k+r} =
n+k-r,T_{n+k+r}=k-n-r \},
$$
since the event $\{ S_{n+k+r} = n+k-r, T_{n+k+r} = k-n-r \}$ implies
that exactly $r$ of $\{ \xi_1 = -1 \}, \cdots, \{ \xi_{n+k+r} = -1
\}$ and $k$ of $\{ \eta_1 = 1 \}, \cdots, \{ \eta_{n+k+r} = 1 \}$
occur, it follows that $ T_{2k-1} \leq 1, T_{2k} \leq 0, T_{2k+j} <
0, j \geq 1 $. Noticing $S_{2k-1} \geq 1$, then we can derive that
$$
\{ S_j \geq 0, S_j \geq T_j, 1 \leq j \leq n+k+r, S_{n+k+r} = n+k-r,
T_{n+k+r} = k-n-r \}
$$
$$
= \{ S_j \geq 0, S_1 \geq T_1, \cdots, S_{2k-2} \geq T_{2k-2},
S_{n+k+r} = n+k-r, T_{n+k+r} = k-n-r \},
$$
that is
$$
\hspace{-2cm} F(n+k+r;n+k-r,n) = 4^{n+k+r} \sum_{s=0}^k P \{ S_j
\geq 0, S_1 \geq T_1, \cdots,
$$
$$
S_{2k-2} \geq 2s+2-2k, T_{2k-2} = 2s+2-2k, S_{n+k+r} = n+k-r,
T_{n+k+r} = k-n-r \}.
$$

\paragraph*{}

Note that the event $\{ T_{2k-2} = 2s+2-2k, T_{n+k+r} = k-n-r \}$ is
equivalent to $\{ T_{2k-2} = 2s+2-2k, \hspace{5pt} exactly
\hspace{5pt} k-s \hspace{5pt} of \hspace{5pt} \{\eta_{2k-1} =1 \},
\cdots, \{\eta_{n+k+r} =1 \} \hspace{5pt} occur \}$, from the
property of i.i.d of $(\xi_1, \eta_1), (\xi_2, \eta_2), \cdots $, we
have

$$
F(n+k+r;n+k-r,n)= 4^{n+k+r} \sum_{s=0}^k \bigg[ {n-k+r+2 \choose
k-s} ( \frac{1}{2}  )^{n-k+r+2} \times
$$
$$
P \{ S_j \geq 0, S_1 \geq T_1, \cdots, S_{2k-2} \geq T_{2k-2} =
2s+2-2k, S_{n+k+r} = n+k-r \} \bigg].
$$

Now we consider the probability above, supposing $S_{2k-2}=2t+2$, to
guarantee $S_{2k-2} \geq T_{2k-2} = 2s+2-2k$ it has to be $-1 \leq t
\leq k-2$ if $s \leq k-1$, and $0 \leq t \leq k-2$ if $s=k$, so that
it follows
$$
P \{ S_j \geq 0, S_1 \geq T_1, \cdots, S_{2k-2} \geq T_{2k-2} =
2s+2-2k, S_{n+k+r} = n+k-r \}
$$
$$
= \sum_t  P \{ S_j \geq 0, S_1 \geq T_1, \cdots, S_{2k-2} = 2t+2,
T_{2k-2} = 2s+2-2k \} \times
$$
$$
P \{ \widetilde{S}_j \geq 0, 1 \leq j \leq n-k+r+2, \widetilde{S}_0
= 2t+2, \widetilde{S}_{n-k+r+2} = n+k-r \},
$$
where $\widetilde{S}_j = 2t + 2 + \xi_{2k-1} + \cdots +
\xi_{2k+j-2}$. From Lemma $1$ the first probability is
$(\frac{1}{4})^{2k-2} F(2k-2;2t+2,k-s+t)$, and by reflection
principle it is easy to obtain the second probability is equal to
$$
(\frac{1}{2})^{n-k+r+2} \left [ {n-k+r+2 \choose n-t} - {n-k+r+2
\choose n+t+3} \right ],
$$
and so on we have
$$
F(n+k+r;n+k-r,n) = \sum_{t=0}^{k-2} \left[ {n-k+r+2 \choose n - t} -
{n-k+r+2 \choose n+t+3} \right] F(2k-2;2t+2,t)
$$
$$
+ \sum_{s=0}^{k-1} \sum_{t=-1}^{k-2}{n-k+r+2 \choose k-s} \left[
{n-k+r+2 \choose n-t} - {n-k+r+2 \choose n+t+3} \right]
F(2k-2;2t+2,k-s+t),
$$
then the formula $(17)$ is proven.

\paragraph*{}

Now we consider the conjecture {\bf C2} of Petkov\v{s}ek and Wilf.
The first formula of $(6.1)$ is immediately obtained from $(10)$. To
verify the other formulas, substituting $r=n+k$ in $(17)$, for $k
\geq 1$, we have
$$
\hspace{-1cm} F(2n + 2k;0,n) = \sum_{t=0}^{k-2} \frac{2t+3}{2n+3}
{2n+3 \choose n+ t+3} F(2k-2;2t+2,t)
$$
$$
+ \sum_{s=0}^{k-1} \sum_{t=0}^{k-1} \frac{2t+1}{2n+3} {2n+2 \choose
s+1} {2n+3 \choose n+t+2} F(2k-2;2t,s+t). \eqno(19)
$$

\vspace{10pt}

(i). Let $k=1$ in (19), note that $F(0;0,0)=1$, it follows that
$$
F(2n+2;0,n) = \frac{2n+2}{n+2} {2n+2 \choose n+1} =
\frac{2^{2n+1}(n+1)(3/2)_n}{(3)_n}.
$$

\vspace{10pt}

(ii). From $(1)-(3)$, it follows that
$$
\begin{array}{l} F(2;0,0)=2,\hspace{3pt} F(2;0,1)=1, \\
F(2;2,0)=1,\hspace{3pt} F(2;2,1)=2,\hspace{3pt} F(2;2,2)=1.
\end{array}
$$

Substituting $k=2$ in $(19)$, we have
$$
F(2n+4;0,n) = \frac{(2n+2)!}{n!
(n+3)!}(8n^2+32n+33)=\frac{4^n(n+1)(8n^2 + 32n +
33)(3/2)_n}{3(4)_n}.
$$

\vspace{10pt}

(iii). From $(1)-(3)$, it follows that
$$
\begin{array}{l} F(4;0,0)=11, \hspace{2pt} F(4;0,1)=8,\hspace{6pt} F(4;0,2)=2, \\
F(4;2,0)=9,\hspace{6pt} F(4;2,1)=17,\hspace{3pt} F(4;2,2)=12,\hspace{3pt} F(4;2,3)=3, \\
F(4;4,0)=1,\hspace{6pt} F(4;4,1)=4,\hspace{8pt}
F(4;4,2)=6,\hspace{8pt} F(4;4,3)=4,\hspace{3pt} F(4;4,4)=1.
\end{array}
$$

Substituting $k=3$ in $(19)$, we have
$$
F(2n+6;0,n) = \frac{(2n+2)!}{3n!(n+4)!} (64n^4 + 672n^3 + 2648n^2 +
4641n + 3060)
$$
$$
=\frac{4^{n-1}(n+1)(64n^4 + 672n^3 + 2648n^2 + 4641n +
3060)(3/2)_n}{9(5)_n}.
$$

Finally, for proving $(6.2)$, it suffices to show that $\frac{n! k!
(n+k+1)!}{(2n+2)!} F(2n+2k;0,n)$ is a polynomial of degree $2k-2$
with all integer coefficients. From the formula $(19)$, it is
immediately obtained because we have the following decomposition
$$
\frac{n! k! (n+k+1)!}{(2n+2)!} F(2n+2k;0,n) = k! \sum_{t=0}^{k-2}
(2t+3) \frac{(n+k+1)!}{(n+t+3)!} \frac{n!}{(n-t)!} F(2k-2;2t+2,t)
$$
$$
+\sum_{s=0}^{k-1} \sum_{t=0}^{k-1} (4t+2) \frac{k!}{(s+1)!}
\frac{(2n+1)!}{(2n-s+1)!} \frac{(n+k+1)!}{(n+t+2)!}
\frac{(n+1)!}{(n-t+1)!} F(2k-2;2t,s+t),
$$
moreover, combining $(10)$ it is not difficult to see $2^{3k-2}
n^{2k-2}$ is the leading term of this polynomial.

\vspace{5pt}

We complete the proof of the conjecture {\bf C2}. \hspace{10pt}
$\Box$

\vspace{15pt}

{\bf Proof of Theorem 5.} Suppose $n \geq 1$. Recalling Lemma 1, we
note that the event $\{ S_{n+2k} = T_{n+2k} = n \}$ implies that
exactly $k$ of $\{ \xi_1 = -1 \}$, $\cdots $, $ \{ \xi_{n+2k} = -1
\}$ occurs and $k$ of $\{ \eta_1 = -1 \}$, $ \cdots $, $\{
\eta_{n+2k} = -1 \}$ occurs, it follows that $T_{2k-1} \geq -1$ and
$S_j \geq 0$ if $j \geq 2k$, therefore
$$
F(n+2k;n,0) = 4^{n+2k} P \{ S_j \geq 0, S_j \geq T_j, 1 \leq j \leq
n+2k, S_{n+2k} = T_{n+2k} = n \}
$$
$$
= \sum_{s=1}^k \sum_{t=0}^s 4^{n+2k} P \{ ( S_i \geq 0, S_i \geq
T_i, 1 \leq i \leq 2k-2, S_{2k-1} = 2s-1, T_{2k-1} = 2t-1 )
$$
$$
\cap  \hspace{3pt} (S_j \geq T_j, 2k \leq j \leq n+2k, S_{n+2k} =
T_{n+2k} = n ) \}
$$
$$
\hspace{-6cm} = \sum_{s=1}^k \sum_{t=0}^s F(2k-1;2s-1,s-t) 4^{n+1} P
\{ A_{s,t} \},
$$
the last equality is from the independent of $\xi, \eta$ and Lemma
1, here $A_{s,t}$ is the intersection event of the following events
$$
\left \{ \begin{array}{l} 2s-1+\sum_{i=0}^j \xi_{2k+i} \geq
2t-1+\sum_{i=0}^j \eta_{2k+i}, 0 \leq j \leq n-1, \\
2s-1+\sum_{i=0}^n \xi_{2k+i} = 2t-1+\sum_{i=0}^n \eta_{2k+i}=n.
\end{array}
\right.
$$

Writing
$$
\widetilde{\xi}_{n+1-i} = \xi_{2k+i}, \hspace{5pt}
\widetilde{\eta}_{n+1-i} = \eta_{2k+i}, \hspace{5pt} 0 \leq i \leq
n,
$$
$$
\widetilde{S}_j = \sum_{i=1}^j \widetilde{\xi}_i, \hspace{5pt}
\widetilde{T}_j = \sum_{i=1}^j \widetilde{\eta}_i, \hspace{5pt} 1
\leq j \leq n+1,
$$
it is obvious that the event $A_{s,t}$ is equivalent to
$$
( \widetilde{T}_j \geq \widetilde{S}_j, 1 \leq j \leq n,
\widetilde{T}_{n+1} = n+1-2t, \widetilde{S}_{n+1} = n+1-2s ),
$$
then Lemma $2$ implies the formula $(18)$.

\paragraph*{}

For the conjecture {\bf C4} of Petkov\v{s}ek and Wilf. The first
formula of $(8.1)$ is immediately obtained from $(9)$. Now we verify
the other formulas.

\vspace{10pt}

(i). Let $k=1$ in (18), note that $F(1;1,0)=F(1;1,1)=1$, it follows
that
$$
F(n+2;n,0) = \frac{n+4}{n+2} {n+2 \choose 2} =
\frac{1}{2}(n+1)(n+4).
$$

\vspace{10pt}

(ii). From $(1)-(3)$, it follows that
$$
\begin{array}{l} F(3;1,0)=5,\hspace{3pt} F(3;1,1)=6, \\
F(3;3,0)=1,\hspace{3pt} F(3;3,1)=3,\hspace{3pt} F(3;3,2)=3.
\end{array}
$$

Substituting $k=2$ in $(18)$, we have
$$
F(n+4;n,0) = \frac{5n+22}{n+2} {n+2 \choose 2} +
\frac{n^2+15n+44}{2(n+2)} {n+2 \choose 3}
$$
$$
= \frac{1}{12}(n+1)(n^3 + 15n^2 + 74n + 132).
$$

\vspace{10pt}

(iii). From $(1)-(3)$, it follows that
$$
\begin{array}{l} F(5;1,0)=37, \hspace{8pt} F(5;1,1)=48, \\
F(5;3,0)=14,\hspace{6pt} F(5;3,1)=36,\hspace{3pt} F(5;3,2)=39, \\
F(5;5,0)=1,\hspace{10pt} F(5;5,1)=5,\hspace{8pt}
F(5;5,2)=10,\hspace{8pt} F(5;5,3)=10.
\end{array}
$$

Substituting $k=3$ in $(18)$, then $F(n+6;n,0)$ is equal to
$$
\frac{37n+170}{n+2} {n+2 \choose 2} + \frac{7n^2+93n+275}{n+2} {n+2
\choose 3} + \frac{n^3+33n^2+272n+660}{6(n+2)} {n+2 \choose 4}
$$
$$
=\frac{1}{144}(n+1)(n^5 + 32n^4 + 407n^3 + 2620n^2 + 8604n + 12240).
$$

Finally, for proving $(8.2)$, it suffices to show that $\frac{k!
(k+1)!}{n+1} F(n+2k;n,0)$ is a polynomial with all integer
coefficients which leading term is $n^{2k-1}$. From the formula
$(18)$, it is clear because it can be decomposed as follows
$$
\sum_{s=1}^k \sum_{t=0}^s (s-t+1) \frac{k!}{t!}
\frac{(k+1)!}{(s+1)!} \frac{(n+2)!}{(n+2-t)!} \frac{n!}{(n+1-s)!}
F(2k-1;2s-1,s-t).
$$

\vspace{5pt}

The proof of the conjecture {\bf C4} is complete. \hspace{10pt}
$\Box$

\section{On the conjectures of Ira Gessel}

\paragraph*{}

It is well known the vicious walker model was introduced by M. E.
Fisher \cite{FI} in $1984$ and received much interest in statistical
physics and in combinatorics, it is related to many important
combinatorial objects, such as partitions, Young tableaux, symmetric
functions...(see [6-8]). The $p$-vicious walker model of length $m$
consists of $p$ lattice paths $W_1, W_2, \cdots, W_p$ in
$\mathbb{Z}^2$ where

\vspace{5pt}

(i). $W_i$ starts from $(0, a_i)$ and ends at $(m, b_i)$ for $1 \leq
i \leq p$,

(ii). all steps are North-East and South-East,

(iii). if $i \neq j$, then $W_i$ and $W_j$ never intersect, in other
words, they do not share any common points.

\vspace{5pt}

In $p$-star configurations, $a_i = 2i-2$ for each $i$ (with no
constraint on $b_i$); In $p$-watermelon configurations, $b_i =
k+2i-2$ for some $k$ \cite{GU}. A vicious walkers configuration with
a wall means it has additional property that none of the paths ever
goes below the $x$-axis. In Section $2$ we have shown that
$F(m;n_1,n_2)$ is the number of pairs of non-crossing lattice paths
$(S,T)$ of length $m$ starting from $(0,0)$ with steps of the form
North-East and South-East, where $S$ running above (may touching)
the $x$-axis and ending at $(m,n_1)$, $T$ running below (may
touching) $S$ ending at $(m,n_1-2n_2)$. If shifting the path $S$ by
$2$ units up, we will obtain a $2$-star configuration, and it is
with half-wall for only $S$ running above (may touching) the
$x$-axis. In particular, Gessel number $F(2n;0,0)$ is the number of
2-watermelons of length $2n$ with half-wall, it is also the number
of pairs of non-crossing Dyck paths and free Dyck paths of length
$2n$.

\paragraph*{}

We note that it is possible to give a human proof of Gessel
conjecture by applying the theory of vicious walkers. From Lemma $1$
we write
$$
F(2n;0,0) = P_1(n) + P_2(n) + P_3(n),
$$
where
$$
\begin{array}{l} P_1(n) = 4^{2n} P \{ S_j \geq 0, T_j \leq 0, 1 \leq j \leq 2n, S_{2n} = T_{2n} =
0 \}, \\ P_2(n) = 4^{2n} P \{ S_j \geq T_j \geq 0, 1 \leq j \leq 2n,
S_{2n} = T_{2n} = 0 \},  \\ P_3(n) = 4^{2n} P \{ S_j \geq 0, S_j
\geq T_j, T \hspace{3pt} crossing \hspace{3pt} x\!\!-\!\!axis,
S_{2n} = T_{2n} = 0 \}.
\end{array}
$$

\vspace{5pt}

Both $P_1(n)$ and $P_2(n)$ are known. It is obvious that
$$
P_1(n) = \frac{{2n \choose n}^2}{n+1}, \eqno(20)
$$
and from the enumeration results of $p$-watermelons with a wall
(\cite{KG}, Theorem 6), it follows that
$$
P_2(n) = \frac{{2n \choose n} {2n+2 \choose n}}{(n+1) {n+3 \choose
n}}, \eqno(21)
$$
the remaining work is to evaluate $P_3(n)$ by using
Lindstr\"{o}m-Gessel-Viennot theorem, it involves summations of
certain binomial coefficients term.

\paragraph*{}

Another possible method is as follows. For $0 \leq k \leq n$, we
write
$$
F_k(n)= 4^{2n} P \{ S_j \geq 0, S_j \geq T_j, 1 \leq j \leq 2n,
S_{2n} = T_{2n} = 0;
$$
$$
\hspace{2cm} exactly \hspace{3pt} 2k \hspace{3pt} of \hspace{3pt}
(\xi_1 = \eta_1), \cdots, (\xi_{2n} = \eta_{2n}) \hspace{3pt} occur
\},  \eqno(22)
$$
then,
$$
F(2n;0,0) = \sum_{k=0}^n F_k(n).
$$

We note that $F_0(n) = \frac{1}{n+1} {2n \choose n} = C_n$, the
$n$-th Catalan number, because $( none \hspace{3pt} of
\\ \hspace{3pt} (\xi_1 = \eta_1), \cdots, (\xi_{2n} = \eta_{2n})
\hspace{3pt} occurs )$ in $(22)$ implies that $T_j = -S_j$ for all
$1 \leq j \leq 2n$. Similarly, $F_n(n)=C_n$. Using words with an
alphabet consisting of four letters, A. Ayyer \cite{AA} also have
obtained the following result:
$$
F_{n-1}(n) = \frac{2n+1}{2} {2n \choose n} - 2^{2n-1}. \eqno(23)
$$

\paragraph*{}

Finally, we mention that there is another conjecture of Gessel
\cite{SL}, denote $a_m$ the number of Gessel walks of length $m$ on
square lattice, starting at origin and staying on points with $x
\geq 0, y \leq x$, it seems to have
$$
a_0=1, \hspace{5pt} a_{2m}= \frac{12m+2}{3m+1} a_{2m-1},
\hspace{5pt} a_{2m+1}= \frac{4m+2}{m+1} a_{2m}. \eqno(24)
$$
As Lemma $1$ follows from the fact that the distribution of $(\xi,
\frac{\xi - \eta}{2} )$ is identical with $(11)$ of $(x,y)$ in
Section $2$, without difficulty we have the following representation
$$
a_m = 4^m P \{ S_j \geq 0, \hspace{5pt} S_j \geq T_j, \hspace{5pt} 1
\leq j \leq m \}
$$
$$ = \sum_{n_1 \equiv m (mod \hspace{2pt} 2) }^m \sum_{n_2 = 0}^{\frac{m+n_1}{2}} F(m;n_1,n_2),
\eqno(25)
$$
This conjecture $(24)$ is the enumeration problem of $2$-star
configurations with arbitrary end points and with half-wall.

\end{document}